\documentclass{article}%
\pdfoutput=1
\usepackage{amsmath}
\usepackage{amsfonts}
\usepackage{amssymb}
\usepackage{graphicx}
\usepackage{graphicx}%
\usepackage{tikz}
\usetikzlibrary{decorations.markings}
\usetikzlibrary{arrows}
\providecommand{\U}[1]{\protect\rule{.1in}{.1in}}
\newtheorem{theorem}{Theorem}

\newtheorem{corollary}[theorem]{Corollary}

\newtheorem{definition}[theorem]{Definition}

\newenvironment{proof}[1][Proof]{\noindent\textbf{#1.} }{\ \rule{0.5em}{0.5em}}
\begin{document}

\title{A note on Dehn colorings and invariant factors}
\author{Derek A. Smith and Lorenzo Traldi\\Lafayette College\\Easton Pennsylvania 18042, United States
\and William Watkins\\California State University Northridge\\Northridge, California 91330, United States}
\date{ }
\maketitle

\begin{abstract}
If $A$ is an abelian group and $\phi$ is an integer, let $A(\phi)$ be the subgroup of $A$ consisting of elements $a \in A$ such that $\phi \cdot a=0$. We prove that if $D$ is a diagram of a classical link $L$ and $0=\phi_0,\phi_1,\dots,\phi_{n-1}$ are the invariant factors of an adjusted Goeritz matrix of $D$, then the group $\mathcal{D}_{A}(D)$ of Dehn colorings of $D$ with values in $A$ is isomorphic to the direct product of $A$ and $A=A(\phi_{0}),A(\phi_1),\dots,A(\phi_{n-1})$. It follows that the Dehn coloring groups of $L$ are isomorphic to those of a connected sum of torus links $T_{(2,\phi_1)} \text{ }\# \text{ } \cdots \text{ } \# \text{ } T_{(2,\phi_{n-1})}$.
\end{abstract} 

\section{Introduction\label{sec:intro}}

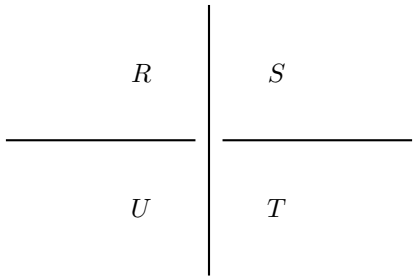
\begin{figure} [bht]
\centering
\begin{tikzpicture} [scale=0.9]
\draw [thick] (-3,0) -- (-0.2,0);
\draw [thick] (0.2,0) -- (3,0);
\draw [thick] (0,2) -- (0,-2);
\node at (-1,1) {$R$};
\node at (1,1) {$S$};
\node at (1,-1) {$T$};
\node at (-1,-1) {$U$};
\end{tikzpicture}
\caption{The regions incident at a crossing.}
\label{figzero}
\end{figure}

If $L=K_{1}\cup\dots\cup K_{\mu}$ is a classical link in $\mathbb{S}^{3}$ then
a \emph{diagram} $D$ of $L$ in the plane is obtained from a regular projection
-- i.e., a projection under which the image of $L$ consists of smooth curves,
with only finitely many singularities, each singularity being a transverse
double point called a \emph{crossing} -- by removing a short segment of the
underpassing on each side of each crossing. Let $\mathcal{R}(D)$ denote the
set of complementary regions of $D$ in the plane. For each crossing $c$ of
$D$, label the incident regions as $R,S,T,U$ in such a way that the pairs of
regions that are neighbors across the overcrossing arc are $\{R,S\}$ and
$\{T,U\}$, as indicated in Figure \ref{figzero}. (If $c$ is not incident on four distinct regions then $R$ or $S$ is an element of $\{T,U\}$.) 

If $A$ is an abelian group then a \emph{Dehn coloring} of $D$ with values in $A$ is a function $d:\mathcal{R}(D)\rightarrow A$ such that $d(R)+d(S)=d(T)+d(U)$ at each crossing. The \emph{Dehn coloring group} of $D$ with values in $A$ is the set $\mathcal{D}_{A}(D)$ of Dehn colorings with values in $A$; it is an abelian group under pointwise addition of functions. Elementary arguments using the Reidemeister moves show that the Dehn coloring groups arising from different diagrams of the same link type are isomorphic. Dehn colorings were introduced by Carter, Silver and Williams \cite{CSW}, and studied further by\ Lamey, Silver and Williams \cite{LSW}. (Our definition differs slightly from \cite{CSW, LSW}, as we do not require a Dehn coloring to be $0$ anywhere.) These authors were expanding on an idea mentioned by Kauffman \cite{K} which, in turn, was suggested by an exercise in the textbook of Crowell and Fox \cite[Exercise VI.6]{CF}. Crowell and Fox defined colorings on the set of arcs of a diagram, rather than the set of regions; these arc colorings are usually called Fox colorings in the literature. We do not discuss Fox colorings in detail, as they are determined by the Dehn colorings \cite{CSW, LSW, Tcol}.

Crowell and Fox introduced link colorings in \cite{CF} as a convenient way to describe certain representations of link groups. Another algebraic interpretation of $\mathcal{D}_{A}(D)$ is provided by the Dehn coloring version of a theorem of Nanyes \cite{N, Tcol}, which implies that $\mathcal{D}_{A}(D)$ is determined up to isomorphism by an adjusted Goeritz matrix of $D$. (See Section 3.) The purpose of the present note is to discuss a structural description of $\mathcal{D}_{A}(D)$, which follows from Nanyes's theorem. The structural description is justified using an argument of Berman \cite{B}, who related the bicycles of a connected graph to invariant factors of the graph's Laplacian matrix. In effect, allowing arbitrary link diagrams involves extending Berman's result to disconnected graphs with signs on their edges; the extra generality does not necessitate any significant change in the argument.

\begin{theorem}
\label{main} Let $\phi_{0},\dots,\phi_{n-1}$ be the invariant factors of an adjusted Goeritz matrix corresponding to one of the two checkerboard shadings of a classical link diagram $D$. Then for every abelian group $A$,
\[
\mathcal{D}_{A}(D)\cong A\times
{\displaystyle\prod\limits_{j=0}^{n-1}}
A(\phi_{j}).
\]
\end{theorem}

Here $A(\phi_j) = \{a \in A \mid \phi_j \cdot a=0\}$. Notice that if $\phi_j=1$ then $A(\phi_j)$ does not make a significant contribution to the direct product of Theorem \ref{main}.

Two corollaries of Theorem~\ref{main} help us understand the amount of information about a link that is reflected in its Dehn coloring groups.

\begin{corollary}
\label{cor}If $L_{1}$ and $L_{2}$ are classical links, then any one of the following conditions implies the other two.

\begin{enumerate}

\item For some diagrams $D_{1}$ and $D_{2}$ of $L_{1}$ and $L_{2}$,
$\mathcal{D}_{A}(D_{1})\cong\mathcal{D}_{A}(D_{2})$ for every abelian group
$A$.

\item For some shadings of diagrams $D_{1}$ and $D_{2}$ of $L_{1}$ and $L_{2}$, the corresponding adjusted Goeritz matrices have the same invariant factors, except that one of the matrices may have extra invariant factors equal to $1$.

\item For some shadings of diagrams $D_{1}$ and $D_{2}$ of $L_{1}$ and $L_{2}$, the abelian groups represented by adjusted Goeritz matrices of $D_1$ and $D_2$ are isomorphic.

\end{enumerate}
\end{corollary}

Corollary \ref{cor} is still valid with \textquotedblleft For some\textquotedblright\ replaced by \textquotedblleft For all\textquotedblright\ in any of the three conditions.

\begin{corollary}
\label{lastcor}
Let $\sigma$ be a shading of a diagram $D$ of a classical link $L$, and let $0=\phi_0,\phi_{1},\dots,\phi_{n-1}$ be the invariant factors of the corresponding adjusted Goeritz matrix. Then there is a connected sum of torus links
\[
T=T_{(2,\phi_1)} \text{ }\# \text{ } \cdots \text{ } \# \text{ } T_{(2,\phi_{n-1})}
\]
such that $\mathcal{D}_{A}(D) \cong  \mathcal{D}_{A}(D_T)$ for every abelian group $A$ and every diagram $D_T$ of $T$.
\end{corollary}

Note that the number of torus links in the connected sum of Corollary \ref{lastcor} is $n-1$, not $n$.

Corollaries \ref{cor} and \ref{lastcor} tell us that if we consider two links with equivalent adjusted Goeritz matrices to be related, then we obtain a rather coarse equivalence relation. This observation contrasts with the fact that if we consider two links with \emph{identical} adjusted Goeritz matrices to be related, then we obtain a much finer equivalence relation; it is the same as the equivalence relation defined by mutation \cite{T}.

Here is an outline of the rest of the note. In Section 2 we recall the theory of invariant factors of matrices over $\mathbb{Z}$. In Section 3 we deduce Theorem \ref{main} and Corollary \ref{cor} from this theory and the theorem of Nanyes \cite{N, Tcol}. Corollary \ref{lastcor} is proven in Section 4, and in Section 5 we mention some relevant references.

Before proceeding we should thank an anonymous reader, whose advice significantly improved our exposition.

\section{Invariant factors and the Smith normal form}

In this section we summarize the theory of the Smith normal form of matrices over $\mathbb{Z}$. The Smith normal form is discussed in many algebra texts, like~\cite[Chapter 3]{J}, \cite[Chapter II]{Ne} and \cite[Chapter 9]{Se}.

\begin{definition}
\label{eld}Suppose $M$ is a matrix with $\kappa$ columns and $\rho$ rows,
whose entries lie in a commutative ring with unity $R$. Then the
\emph{elementary ideals} $E_j(M)$ are ideals of $R$, indexed by $j \geq 0 \in \mathbb{N}$. If $\kappa \leq j$, then $E_{j}(M)$ is the ideal $(1)$, i.e., the entire ring $R$. If $\max\{0,\kappa-\rho\} \leq j \leq \kappa-1$, then $E_{j}(M)$ is the ideal
generated by the determinants of the $(\kappa-j)\times(\kappa-j)$ submatrices
of $M$. If $0\leq j<\max\{0,\kappa-\rho\}$, then $E_{j}(M)=\{0\}$.
\end{definition}

\begin{definition}
\label{eldiv}
If $R=\mathbb{Z}$ in Definition \ref{eld}, then for each index
$j\geq0$, let $\delta_{j}(M)$ denote the greatest common divisor of the ideal
$E_{j}(M)$. The integers
\[
\phi_{0}(M)=\frac{\delta_{0}(M)}{\delta_{1}(M)}\text{, }\phi_{1}%
(M)=\frac{\delta_{1}(M)}{\delta_{2}(M)}\text{, }\dots\text{, }\phi_{\kappa
-1}(M)=\frac{\delta_{\kappa-1}(M)}{\delta_{\kappa}(M)}%
\]
are the \emph{invariant factors} of $M$. (We use
the notation $\frac{0}{0}=0$.)
\end{definition}

It is not immediately apparent from Definition \ref{eldiv}, but it turns out that the invariant factors form a sequence of divisors: $\phi_{j}(M) \mid \phi_{j-1}(M)$ $\forall j$. Notice that Definition \ref{eldiv} allows invariant factors equal to $0$. Some references, like \cite[Chapter 3]{J} and \cite[Chapter II]{Ne}, do not allow such factors; instead they simply mention that the number of nonzero invariant factors is the rank of $M$ over $\mathbb{Q}$.

In the situation of Definition \ref{eldiv}, the \emph{Smith normal form}
$\Phi(M)$ is the $\rho\times\kappa$ matrix obtained from the diagonal matrix%
\begin{equation*}%
\begin{pmatrix}
\phi_{0}(M) & 0 & 0 & 0\\
0 & \phi_{1}(M) & 0 & 0\\
0 & 0 & \ddots & 0\\
0 & 0 & 0 & \phi_{\kappa-1}(M)
\end{pmatrix}
\end{equation*}
by adjoining $\rho-\kappa$ rows of zeroes if $\kappa<\rho$, and removing
$\kappa-\rho$ rows of zeroes if $\kappa>\rho$.

\begin{theorem}
\label{smith}The Smith normal form over $\mathbb{Z}$ has the following properties:
\end{theorem}

\begin{enumerate}
\item $\Phi(M)$ is equivalent to $M$, i.e., there are unimodular matrices
$U_{1},U_{2}$ such that $U_{1}MU_{2}=\Phi(M)$.

\item $\Phi(M)$ is the unique matrix equivalent to $M$ whose off-diagonal
entries are all $0$ and whose diagonal entries $M_{jj}$ are non-negative
integers such that $M_{jj}\mid M_{(j-1)(j-1)}$ $\forall j$.

\item If $M^{\prime}$ is another $\rho\times\kappa$ matrix of integers, then any one of the following statements holds only if all four hold.

\begin{enumerate}
\item $M$ and $M^{\prime}$ are equivalent.

\item $\Phi(M)=\Phi(M^{\prime})$.

\item The homomorphisms $f_M,f_{M'}:\mathbb{Z}^{\rho}\rightarrow\mathbb{Z}^{\kappa}$
represented by $M$ and $M^{\prime}$ have isomorphic cokernels, i.e.,  $\mathbb{Z}^{\kappa}/f_M(\mathbb{Z}^{\rho}) \cong \mathbb{Z}^{\kappa}/f_{M'}(\mathbb{Z}^{\rho})$. 

\item For every abelian group $A$, the homomorphisms $A^{\rho}\rightarrow
A^{\kappa}$ defined by $M$ and $M^{\prime}$ have isomorphic kernels.
\end{enumerate}
\end{enumerate}

\begin{proof}
The only part of the theorem that is not mentioned in standard texts like \cite[Chapter 3]{J} or \cite[Chapter II]{Ne} is the equivalence of condition
(d) with the other conditions of property 3. We proceed to explain this equivalence. For (a) $\implies$ (d), notice that if there are unimodular matrices
$U_{1},U_{2}$ such that $U_{1}MU_{2}=M^{\prime}$ then $U_{1}$ and $U_{2}$
represent automorphisms of $A^{\rho}$ and $A^{\kappa}$, so the homomorphisms
$A^{\rho}\rightarrow A^{\kappa}$ defined by $M$ and $M^{\prime}$ have
isomorphic kernels. For (d) $\implies$ (b), notice that if $A$ is an abelian
group then the homomorphism $A^{\rho}\to A^{\kappa}$ represented by
$M$ differs from the homomorphism $A^{\rho}\to A^{\kappa}$ represented by $\Phi(M)$ only through composition with automorphisms of the domain and codomain; so the homomorphisms represented by $M$ and $\Phi(M)$ certainly have isomorphic kernels. Requiring that the nonzero entries of $\Phi(M)$ all lie on the diagonal, and the diagonal entries of $\Phi(M)$ form a sequence of divisors, implies that the diagonal entries of $\Phi(M)$ are completely determined by the kernels of the homomorphisms $A^{\rho} \to A^{\kappa}$ represented by $\Phi(M)$ for $A=\mathbb{Q}$ and $A=\mathbb{Z}_{m}$, $m>1\in\mathbb{N}$. The same holds for $M^{\prime}$, so (d) implies (b).
\end{proof}

\section{Theorem \ref{main} and Corollary \ref{cor}}

It is well known that the complementary regions of a link diagram $D$ can be separated into two classes -- \emph{shaded} and \emph{unshaded} -- in such a way that regions whose boundaries share an arc are shaded differently. We think of the shaded/unshaded designations as arbitrary, so $D$ has two different shadings, which are opposites of each other. We do not prefer one shading or the other. A shading $\sigma$ of $D$ has two associated \emph{checkerboard graphs}, which we denote $\Gamma_{s}(D,\sigma)$ and $\Gamma_{u}(D,\sigma)$. The former has a vertex for each shaded region, and the latter has a vertex for each unshaded region. Each graph has an edge for each crossing, incident on the vertex or vertices corresponding to region(s) which appear at that crossing.

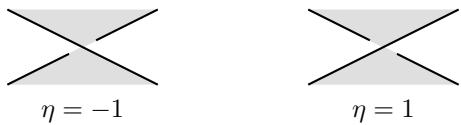
\begin{figure} [bht]
\centering
\begin{tikzpicture} 
\draw [lightgray!50, fill=lightgray!50] (-2.5,0.5) -- (-1.5,0) -- (-.5,0.5);
\draw [lightgray!50, fill=lightgray!50] (-2.5,-0.5) -- (-1.5,0) -- (-1/2,-0.5);
\draw [lightgray!50, fill=lightgray!50] (3.5,0.5) -- (5/2,0) -- (1.5,0.5);
\draw [lightgray!50, fill=lightgray!50] (7/2,-0.5) -- (5/2,0) -- (3/2,-0.5);
\draw [thick] (-2.5,0.5) -- (-1/2,-0.5);
\draw [thick] (-2.5,-0.5) -- (-1.68,-.09);
\draw [thick] (-1.32,.09) -- (-1/2,0.5);
\draw [thick] (3.5,0.5) -- (1.5,-0.5);
\draw [thick] (3.5,-0.5) -- (2.68,-.09);
\draw [thick] (2.32,.09) -- (1.5,0.5);
\node at (-3/2,-0.85) {$\eta=-1$};
\node at (5/2,-0.85) {$\eta=1$};
\end{tikzpicture}
\caption{The Goeritz index.}
\label{figone}
\end{figure}

Suppose $D$ is a link diagram with a shading $\sigma$, and $U_{1},\dots,U_{n}$
are the unshaded complementary regions with respect to $\sigma$. For $i\neq j$
let $C_{ij}$ be the set of crossings that are incident on both $U_{i}$ and
$U_{j}$. For each crossing, define the \emph{Goeritz index} $\eta$ as in
Figure \ref{figone}.

Then the Goeritz matrix $G(D,\sigma)$ is the symmetric $n\times n$ matrix whose entries are given by%
\[
G(D,\sigma)_{ij}=%
\begin{cases}
-\sum\limits_{c\in C_{ij}}\eta(c)\text{,} & \text{if }i\neq j\\
-\sum\limits_{k\neq i}G(D,\sigma)_{ik}\text{,} & \text{if }i=j
\end{cases} \text{.}
\]

The Goeritz matrix was introduced more than 80 years ago \cite{G}. That was a
time when knot theorists did not pay much attention to links with $\mu>1$, so
the fact that $G(D,\sigma)$ is insensitive to simple link operations (like
inserting a crossing-free circle into an unshaded region of a link diagram)
was not a cause for concern. To incorporate more information about links, it
is convenient to adjust the definition of the matrix as follows \cite{Told, T}.

\begin{definition}
\label{adj}Let $D$ be a link diagram with a shading $\sigma$, and let
$\beta_{s}(D,\sigma)$ denote the number of connected components of the shaded
checkerboard graph $\Gamma_{s}(D,\sigma)$. If $\beta_{s}(D,\sigma)>1$ then the
\emph{adjusted Goeritz matrix} of $D$ with respect to $\sigma$ is%
\[
G^{adj}(D,\sigma)=%
\begin{pmatrix}
G(D,\sigma) & 0\\
0 & 0
\end{pmatrix}
\]
where the $0$ in the lower right denotes a $(\beta_{s}(D,\sigma)-1)\times
(\beta_{s}(D,\sigma)-1)$ zero matrix. If $\beta_{s}(D,\sigma)=1$ then
$G^{adj}(D,\sigma)=G(D,\sigma)$.
\end{definition}

We are now ready to prove Theorem \ref{main}. Let $A$ be an abelian group, let $\sigma$ be a shading of a link diagram $D$, and let $f_G,f_{\Phi}$ be  the homomorphisms
\[A^{n+\beta_{s}(D,\sigma)-1} \to A^{n+\beta_{s}(D,\sigma)-1}\] 
represented by $G^{adj}(D,\sigma)$ and $\Phi(G^{adj}(D,\sigma))$, respectively. The Dehn coloring version of the theorem of Nanyes~\cite{N, Tcol} is equivalent to the assertion that $\mathcal{D}_{A}(D)$ is isomorphic to the direct product $A \times \ker f_G$.  Theorem \ref{smith} tells us that $f_G$ and $f_{\Phi}$ are related to each other through composition with automorphisms of $A^{n+\beta_{s}(D,\sigma)-1}$, so the kernels of $f_G$ and $f_{\Phi}$ are isomorphic. It follows that $\mathcal{D}_{A}(D)$ is isomorphic to $A \times \ker f_{\Phi}$. This isomorphism completes the proof of Theorem \ref{main}. 

Corollary \ref{cor} follows from Theorems \ref{main} and \ref{smith}. The ``for all'' versions of the conditions of Corollary \ref{cor} follow from the fact that if $\sigma,\sigma'$ are shadings of diagrams $D,D'$ of the same link $L$, then the resulting adjusted Goeritz matrices are related to each other through operations whose only effect on the Smith normal form is to change the number of invariant factors equal to $1$. (See \cite{P, Told} for detailed explanations of this fact.)

Before proceeding we should mention that there is also a version of Theorem \ref{main} for Fox colorings; the proof is much the same. Using the definitions and notation of \cite{Tcol}, this version of Theorem \ref{main} is
\[
\mathcal{F}_{A}(D)\cong 
{\displaystyle\prod\limits_{j=0}^{n-1}}
A(\phi_{j}).
\]

\section{Corollary \ref{lastcor}}

Let $L$ be a link with an adjusted Goeritz matrix whose invariant factors are $\phi_0,\dots,\phi_{n-1}$. The columns and rows of $G^{adj}(D,\sigma)$ sum to $0$, so $\phi_0=0$. To build the diagram $D_T$ of Corollary \ref{lastcor}, begin with a simple circle, which bounds a shaded region. Use connected sums to attach torus links $T_{(2,\phi_1)},\dots, T_{(2,\phi_{n-1})}$ inside the circle, one after the other. The resulting diagram $D_T$ and shading $\sigma_T$ have these properties: $\beta_s(D_T,\sigma_T)=1$, there are $n$ unshaded regions, $\phi_j=G^{adj}(D_T,\sigma_T)_{jj}= - G^{adj}(D_T,\sigma_T)_{nj}= - G^{adj}(D_T,\sigma_T)_{jn}$ for $1 \leq j \leq n-1$, $G^{adj}(D_T,\sigma_T)_{nn} = \phi_1+\dots+\phi_{n-1}$, and all other entries of $G^{adj}(D_T,\sigma_T)$ are $0$. 

The proof of Corollary \ref{lastcor} is completed by verifying that the invariant factors of $G^{adj}(D_T,\sigma_T)$ are $0=\phi_0,\phi_1,\dots,\phi_{n-1}$. To see why this is so, add rows $1$ through $n-1$ of $G^{adj}(D_T,\sigma_T)$ to the last row, and then add columns $1$ through $n-1$ to the last column. The resulting matrix is equivalent to $G^{adj}(D_T,\sigma_T)$. It has $\phi_1,\dots,\phi_{n-1},0$ on the diagonal, and all other entries equal to $0$.

\begin{figure} 
\centering
\begin{tikzpicture} [scale=0.38]
\draw [lightgray!50, fill=lightgray!50] (-27,4) -- (-27,0) -- (-26,0) -- (-25.33,0.67) -- (-25.33,1.33) -- (-26,2) -- (-26,4);
\draw [lightgray!50, fill=lightgray!50] (-26,0.67) -- (-26,1.33) -- (-24.67,1.33) -- (-24,2) -- (-23,2) -- (-23,0) -- (-24,0) -- (-24.67,0.67);
\draw [lightgray!50, fill=lightgray!50] (-24,4) -- (-20,4) -- (-20,2) -- (-19,1) -- (-20,0) -- (-24,0);
\draw [lightgray!50, fill=lightgray!50] (-27,4) -- (3,4) -- (3,3) -- (-27,3);
\draw [lightgray!50, fill=lightgray!50] (3,4) -- (3,0) -- (2,0) -- (1,1) -- (2,2) -- (2,4);
\node at (-25,2.2) {1};
\draw [thick] (-27,4) -- (-27,0);
\draw [thick] (-26,0) -- (-27,0);
\draw [thick] (-20,0) -- (-24,0);
\draw [thick] (-24,0) -- (-24.67,0.67);
\draw [thick] (-24.67,0.67) -- (-25.33,0.67);
\draw [thick] (-25.33,0.67) -- (-26,0);
\draw [thick] (-20,0) -- (-24,0);
\draw [thick] (-24,2) -- (-24.67,1.33);
\draw [thick] (-24.67,1.33) -- (-25.33,1.33);
\draw [thick] (-25.33,1.33) -- (-26,2);
\draw [thick] (-26,3) -- (-26,2);
\draw [thick] (-26,3) -- (-24,3);
\draw [thick] (-24,2) -- (-24,3);
\node at (-17,2.2) {2};
\draw [lightgray!50, fill=lightgray!50] (-19,1) -- (-18,0) --(-17,1) -- (-18,2);
\draw [lightgray!50, fill=lightgray!50] (-17,1) -- (-16,0) --(-15,1) -- (-16,2);
\draw [lightgray!50, fill=lightgray!50] (-15,1) -- (-14,0) --(-10,0) -- (-9,1) -- (-10,2) -- (-10,4) -- (-14,4) -- (-14,2);
\draw [thick] (-20,2) -- (-18,0);
\draw [thick] (-20,0) -- (-19.2,0.8);
\draw [thick] (-18.8,1.2) -- (-18,2);
\draw [thick] (-18,2) -- (-16,0);
\draw [thick] (-18,0) -- (-17.2,0.8);
\draw [thick] (-16.8,1.2) -- (-16,2);
\draw [thick] (-16,2) -- (-14,0);
\draw [thick] (-16,0) -- (-15.2,0.8);
\draw [thick] (-14.8,1.2) -- (-14,2);
\draw [thick] (-14,3) -- (-14,2);
\draw [thick] (-14,3) -- (-20,3);
\draw [thick] (-20,2) -- (-20,3);
\node at (-7,2.2) {3};
\draw [lightgray!50, fill=lightgray!50] (-9,1) -- (-8,0) --(-7,1) -- (-8,2);
\draw [lightgray!50, fill=lightgray!50] (-5,1) -- (-6,0) --(-7,1) -- (-6,2);
\draw [lightgray!50, fill=lightgray!50] (1,1) -- (0,0) --(-4,0) -- (-5,1) -- (-4,2) -- (-4,4) -- (0,4) -- (0,2);
\draw [thick] (-14,0) -- (-10,0);
\draw [thick] (-10,2) -- (-8,0);
\draw [thick] (-10,0) -- (-9.2,0.8);
\draw [thick] (-8.8,1.2) -- (-8,2);
\draw [thick] (-8,2) -- (-6,0);
\draw [thick] (-8,0) -- (-7.2,0.8);
\draw [thick] (-6,2) -- (-4,0);
\draw [thick] (-6,0) -- (-5.2,0.8);
\draw [thick] (-6.8,1.2) -- (-6,2);
\draw [thick] (-4.8,1.2) -- (-4,2);
\draw [thick] (-4,3) -- (-4,2);
\draw [thick] (-4,3) -- (-10,3);
\draw [thick] (-10,2) -- (-10,3);
\node at (1,2.2) {4};
\node at (4,2.2) {5};
\draw [thick] (0,0) -- (-4,0);
\draw [thick] (0,2) -- (2,0);
\draw [thick] (0,0) -- (0.8,0.8);
\draw [thick] (2,2) -- (1.2,1.2);
\draw [thick] (2,2) -- (2,3);
\draw [thick] (0,3) -- (2,3);
\draw [thick] (0,3) -- (0,2);
\draw [thick] (3,0) -- (2,0);
\draw [thick] (-27,4) -- (3,4);
\draw [thick] (3,0) -- (3,4);
\end{tikzpicture}
\caption{A connected sum of torus links, $T_{(2,0)} \text{ }\# \text{ }  T_{(2,3)} \text{ }\# \text{ }  T_{(2,3)} \text{ }\# \text{ }  T_{(2,1)}$. The numbers in the figure indicate row and column indices in the Goeritz matrix mentioned after the proof of Corollary \ref{lastcor}.}
\label{figtwo}
\end{figure}
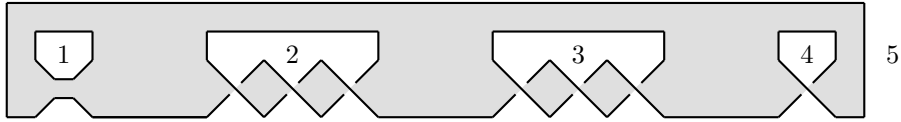

For example, to realize the invariant factors $\phi_0=\phi_1=0$, $\phi_2=\phi_3=3$, and $\phi_4=1$ we construct the diagram $D_T$ pictured in Figure \ref{figtwo}. This diagram has
\[
G^{adj}(D_T,\sigma_T)=
\begin{pmatrix}
0 & 0 & 0 & 0 & 0\\
0 & 3 & 0 & 0 & -3\\
0 & 0 & 3 & 0 & -3\\ 
0 & 0 & 0 & 1 & -1\\
0 & -3 & -3 & -1 & 7
\end{pmatrix} \text{,}
\]
with $\sigma_T$ the indicated shading.

Notice that if $\phi_j=1$ then the corresponding trivial torus link $T_{(2,1)}$ does not make a significant contribution to the connected sum of Corollary \ref{lastcor}. This corresponds to the fact that $A(\phi_j)=A(1)=\{0\}$ does not make a significant contribution to the direct product of Theorem \ref{main}.

\section{Connections with other invariants}

The oldest link invariant connected with the Dehn coloring groups is the abelian group represented by the (adjusted) Goeritz matrix. Seifert \cite{S} proved that for a knot, this group is isomorphic to the direct sum of $\mathbb{Z}$ and the first homology group of the cyclic double cover of $\mathbb{S}^{3}$ branched along the knot; Seifert's result was extended to links by Kyle \cite{Ky}. For more recent accounts see \cite{GL} and \cite[Chapter 9]{L}.

The invariant factors of (adjusted) Goeritz matrices are also connected to other classical link invariants, the elementary ideals. If $L=K_{1} \cup \dots \cup K_{\mu}$ then the group $G$ of $L$ has the property that the abelianization
$G/G^{\prime}$ is free abelian of rank $\mu$. A basis for this free abelian
group consists of meridians $t_{1},\dots,t_{\mu}$ of the components of $L$,
and the integral group ring $\mathbb{Z}[G/G^{\prime}]$ is naturally isomorphic
to the ring $\mathbb{Z}[t_{1}^{\pm1},\dots,t_{\mu}^{\pm1}]$ of Laurent
polynomials in $t_{1},\dots,t_{\mu}$ (with integer coefficients). The
elementary ideals $E_{j}(L)$ are ideals of $\mathbb{Z}[t_{1}^{\pm1},\dots,t_{\mu}^{\pm1}]$ obtained from Definition \ref{eld}, using an Alexander matrix\ of $L$ for $M$. Alexander matrices are discussed in standard references like \cite{CF} and \cite{F}.

Let $\nu:\mathbb{Z}[t_{1}^{\pm1},\dots,t_{\mu}^{\pm1}]\rightarrow\mathbb{Z}$
be the \emph{negative augmentation map}, i.e., the homomorphism of rings with unity given by
$\nu(t_{i}^{\pm1})=-1$ $\forall i$. Then the invariant factors of $G^{adj}(D,\sigma)$ are determined by the images of the elementary ideals of $L$ under $\nu$, because $\nu(E_{j}(L))=E_{j}(G^{adj}(D,\sigma))$ $\forall j$. This fact is implicit in classical discussions of these ideas, like \cite{F}, but we have not found an explicit statement older than \cite{Told}. 

The invariant factors of Goeritz matrices have been studied recently by several authors, including Ikeda and Sugimoto \cite{IS} and Tanaka \cite{Ta}. These authors refer to the sequence of invariant factors (equivalently, the matrix $\Phi(G^{adj}(D,\sigma))$) as the \textquotedblleft
Goeritz invariant\textquotedblright\ of a link, and credit Kawauchi for the idea of studying it as a link invariant. Corollary \ref{cor} tells us that two links have the same Goeritz invariant if and only if the Dehn coloring groups of one link are isomorphic to the Dehn coloring groups of the other link.

\end{document}